\documentclass[12pt]{article}
\usepackage{amssymb}
\usepackage{amsfonts}
\usepackage{amsthm}
\usepackage{amsmath}
\newtheorem{theorem}{Theorem}[section]
\numberwithin{equation}{section}

\begin{document}
\title{Variations of Stieltjes-Wigert and $q$-Laguerre polynomials  and their recurrence coefficients}
\author{Lies Boelen\footnote{current address: Department
of Medicine, Imperial College London, London, United Kingdom}\ \ and  Walter Van Assche\footnote{Department of Mathematics, KU Leuven, Belgium} \\ KU Leuven}
\date{\textit{Dedicated to Dick Askey at the occasion of his 80th birthday}}
\maketitle

\begin{abstract}
We look at some extensions of the Stieltjes-Wigert weight functions. 
First we replace the variable $x$ by $x^2$ in a family of weight functions given by Askey in 1989 and we show that the recurrence coefficients of the corresponding orthogonal polynomials can be expressed in terms of a solution of the $q$-discrete Painlev\'e III equation $q$-$\textrm{P}_{\textrm{III}}$. Next we consider the $q$-Laguerre or generalized Stieltjes-Wigert weight functions with a quadratic transformation and derive recursive equations for the recurrence coefficients of the orthogonal polynomials. These turn out to be related to the $q$-discrete Painlev\'e V equation
$q$-$\textrm{P}_{\textrm{V}}$. Finally we also consider the little $q$-Laguerre weight with a quadratic transformation and show that the recurrence coefficients of the orthogonal polynomials are again related to
$q$-$\textrm{P}_{\textrm{V}}$.
\end{abstract}

\section{Introduction}

In this paper we will investigate semiclassical variations of the Stieltjes-Wigert polynomials, the $q$-Laguerre polynomials and
the little $q$-Laguerre polynomials. We are mainly interested in the recurrence coefficients of the monic polynomials in the three term
recurrence relation
\begin{equation}  \label{TTR}
   xP_n(x) = P_{n+1}(x) + b_n P_n(x) + a_n^2 P_{n-1}(x), \qquad n \geq 0, 
\end{equation}
where $P_0=1$ and $P_{-1}=0$. There are two reasons for combining these three families in one paper. First of all, the Stieltjes-Wigert 
polynomials are really just a special instance of the $q$-Laguerre polynomials and there is also a connection between the $q$-Laguerre polynomials
and the little $q$-Laguerre polynomials (changing $q$ to $1/q$). The second reason is that we use the same technique to prove the results for these three families of
recurrence coefficients, which results in very similar computations. 

In the first three subsections we will briefly recall the Stieltjes-Wigert, $q$-Laguerre and little $q$-Laguerre polynomials. In \S \ref{section:main}
we will introduce the semiclassical variation of these three families and we will formulate our main results.

\subsection{Stieltjes-Wigert polynomials}
Stieltjes introduced in
\cite[\S 56]{Stieltjes} the family of positive weight functions 
\[ \frac{1}{\sqrt{\pi}} x^{-\log x}[1+\lambda \sin (2\pi\log x)], \qquad x \in [0,\infty) ,\quad -1 \leq \lambda \leq 1, \] 
and showed that the moments exist and are independent of $\lambda$:
\[   \frac{1}{\sqrt{\pi}}\int_0^\infty x^n x^{-\log x}[1+\lambda \sin (2\pi\log x)]\, dx = e^{(n+1)^2/4}.  \]
We denote the weight function with $\lambda=0$ by
$$
w_1(x) = \frac{1}{\sqrt{\pi}} \exp(-\log^2(x)), \qquad x \in [0,\infty),
$$
and Stieltjes' main contribution is that this weight function is not uniquely determined by its moments, since there
are other positive measures on $[0,\infty)$ with the same moments: the measure with weight function $w_1$ corresponds to 
an \textit{indeterminate moment problem}. Wigert \cite{Wigert} introduced the more general weight function
\begin{equation}\label{sw-weight}
w_k(x) = \frac{k}{\sqrt{\pi}} \exp(-k^2\log^2(x)), \qquad  x \in [0,\infty),\,\, k>0,
\end{equation}
for which the moments are given by
$$
\mu_n = \int_0^\infty x^n w_k(x)\, dx = q^{-\frac{(n+1)^2}{2}}, 
$$
where $q=\exp(-1/(2k^2))$. This weight function also induces a positive measure on $[0,\infty)$ with an 
indeterminate Stieltjes moment problem (\cite[Section 2.6]{Chihara}).
Another measure with the same moments is given by Askey \cite{Askey}: he
considers weight functions of the form
\begin{equation}\label{askey-weight}
\tilde{w}_\alpha(x) = \frac{x^\alpha}{(-x;q)_\infty(-q/x;q)_\infty},
\qquad x >0,\ 0 < q < 1.
\end{equation}
Here $(x;q)_\infty$ is the $q$-Pochhammer symbol
\[  (x;q)_\infty = \prod_{k=0}^\infty(1-xq^k). \]
The moments $m_n$ associated with this weight function satisfy
$$
\frac{m_n}{m_0}= q^{-n(\alpha+1)-n(n-1)/2}.
$$
This can be compared with the moments of the Stieltjes-Wigert weight $\mu_n/\mu_0=q^{-n^2/2-n}$, and hence for
$\alpha=1/2$ Askey's weight function has the same moments as the Stieltjes-Wigert weight\footnote{in \cite{Askey} there is an unfortunate misprint 
in (4.14) and (4.15) which should be $d(n) = q^{-n^2/2-n-1/2}$ and $d(n)/d(0) = q^{-n^2/2-n}$. Then $\gamma=3/2$ and in (4.16) the $x^{-5/2}$ should 
be $x^{1/2}$.}. In particular this means that
the orthogonal polynomials for the weight function (\ref{askey-weight}) with $\alpha=1/2$ coincide with the Stieltjes-Wigert polynomials, which are
orthogonal with respect to $w_k$. Observe that Askey's weight function has a factor $x^\alpha$ so that (\ref{askey-weight}) 
is in fact a one-parameter family of weight functions with explicitly known moments. This is implicitly also the case for the Stieltjes-Wigert
weight (\ref{sw-weight}), which has the scaling property
\[    w_k(q^\beta x) = q^{\beta^2/2} x^\beta w_k(x), \]
where we used $\log q = -1/(2k^2)$. Chihara \cite{Chihara2} also gives a number of discrete measures with the same moments
as the Stieltjes-Wigert weight. A detailed study of various solutions of the moment problem for the Stieltjes-Wigert weight,
including N-extremal solutions, can be found in the work of Christiansen \cite{Jacob1}.
The orthonormal Stieltjes-Wigert polynomials are given by
\[  P_n(x) = (-1)^n \frac{q^{(2n+1)/4}}{\sqrt{(q;q)_n}}  \sum_{k=0}^n \frac{(q;q)_n}{(q;q)_k(q;q)_{n-k}} q^{k^2} (-q^{1/2}x)^k, \]
\cite[Eq. (11)]{Wigert}, \cite[Eq. (2.3) on p.~173]{Chihara}, \cite[Eq. (2.7.4) on p.~33]{Szego}.
In terms of  basic hypergeometric functions \cite{GasperRahman}, the monic Stieltjes-Wigert polynomials are
\begin{equation}  \label{Pn}
P_n(x) = (-1)^n q^{-n^2-\frac{n}{2}} {}_1\phi_1 \left(\begin{array}{c} q^{-n} \\ 0 \end{array} ;q , -q^{n+3/2}x \right).
\end{equation}
Note that the polynomials in \cite[Eq. (3.27.1) on p.~116]{KoekoekSwarttouw} \cite[Eq. (14.27.1) on p.~544]{KoekoekLeskySwarttouw}
\[    S_n(x) = \frac{1}{(q;q)_n}{}_1\phi_1 \left( \begin{array}{c} q^{-n} \\ 0 \end{array} ; q , -q^{n+1}x \right) \]
are orthogonal for the weight function $x^{-1/2} w_k(x)$ (or Askey's weight function $\tilde{w}_0(x)$)\footnote{The remark in \cite[p.~117]{KoekoekSwarttouw} and \cite[p.~546--547]{KoekoekLeskySwarttouw} is missing this factor $x^{-1/2}$ in the weight function.}.
The monic Stieltjes-Wigert polynomials satisfy the three-term recurrence relation \eqref{TTR}
with coefficients $b_n = q^{-2n-3/2}\left(1+q-q^{n+1}\right)$ and $a_n^2 = q^{-4n}(1-q^n)$. The monic orthogonal polynomials
for the weight $w_k(q^\beta x)=q^{\beta^2/2} x^\beta w_k(x)$ are $Q_n(x) = q^{-n\beta} P_n(q^\beta x)$ and satisfy
\[   xQ_n(x) = Q_{n+1}(x) + b_nq^{-\beta} Q_n(x) + a_n^2q^{-2\beta} Q_{n-1}(x), \]
with $a_n^2$ and $b_n$ the recurrence coefficients of the Stieltjes-Wigert polynomials $P_n$. 
The case $\beta=-1$ gives the log-normal distribution, which is well known in statistics (the distribution of $Y= e^X$, where
$X$ has a standard normal distribution). The case $\beta=\alpha-1/2$ gives a weight function with the same moments
as Askey's weight function $\tilde{w}_\alpha$ in \eqref{askey-weight}. 

\subsection{$q$-Laguerre polynomials}  \label{sec:1.2}
Chihara \cite[pp.~173--174]{Chihara} \cite{Chihara1} \cite{Chihara3} introduced a generalization of the
Stieltjes-Wigert weight function:
\begin{equation}\label{chihara-weight}
\rho_0(x)=  \bigl(\frac{-p}{\sqrt{q}x};q\bigr)_\infty w_k(x),
\qquad x \in [0,\infty),\,\,\, p \in [0,1).
\end{equation}
The monic orthogonal polynomials with respect to \eqref{chihara-weight} are
\[    S_n(x;p,q) = (-1)^n q^{-n(2n+1)/2} (p;q)_n \sum_{k=0}^n \frac{(q;q)_n}{(q;q)_k(q;q)_{n-k}} q^{k^2} \frac{(-q^{1/2}x)^k}{(p;q)_k}, \]
or, in basic hypergeometric form,
\begin{equation} \label{Pn-p}
   S_n(x;p,q) = (-1)^n q^{-n(2n+1)/2} (p;q)_n\ {}_2\phi_1 \left( \begin{array}{c} q^{-n} \\ p \end{array} ; q, -q^{n+\frac32}x \right), 
\end{equation}
which, for $p=0$, reduces to \eqref{Pn}. These orthogonal polynomials are nowadays known as $q$-Laguerre polynomials
\[   L_n^{(\alpha)}(x;q) = \frac{(q^{\alpha+1};q)_n}{(q;q)_n} 
  {}_1\phi_1 \left( \begin{array}{c} q^{-n} \\ q^{\alpha+1} \end{array} ; q, -q^{n+\alpha+1} x \right), \]
\cite[Eq. (3.21.1) on p.~108]{KoekoekSwarttouw}, \cite[Eq. (14.21.1) on p.~522]{KoekoekLeskySwarttouw} and were also investigated by Moak \cite{Moak}.
Indeed, if $p=q^{\alpha+1}$ then $S_n(x;p,q) = C_n L_n^{(\alpha)}(xq^{3/2}/p;q)$, where $C_n$ is a normalizing constant.
The $q$-Laguerre polynomials were also investigated by Askey \cite{Askey2}, Ismail and Rahman \cite{IsmailRahman} and Christiansen \cite{Jacob2}.
Chihara found the recurrence coefficients for the monic orthogonal
polynomials $S_n(x;p,q)$ with respect to the weight \eqref{chihara-weight}:
\begin{align*}
b_n &= q^{-n-3/2}\left(-p-q+(1+q)q^{-n} \right), \\
a_n^2 &= q^{-4n}\left( 1-q^n\right)\left( 1-pq^{n-1}\right),
\end{align*}
which for $p=0$ also reduces to the Stieltjes-Wigert case.
The moment problem for these $q$-Laguerre polynomials is also indeterminate and in fact one has the orthogonality relation
\[  \int_0^\infty L_n^{(\alpha)}(x;q)L_m^{(\alpha)}(x;q) \frac{x^\alpha}{(-x;q)_\infty}\, dx = 0, \qquad n \neq m, \]
which uses a somewhat simpler weight function than \eqref{chihara-weight}.
Christiansen investigated the various solutions of the moment problem for $q$-Laguerre polynomials and one of the weights
he found is 
\[   v(x) =  \frac{(-q/x;q)_\infty}{(-q^{\alpha+1-c}x)_\infty (-q^{-\alpha+c}/x;q)_\infty} x^{c-1}  \]
\cite[Eq. (3.0.6) on p.~9]{Jacob2}.
Indeed $v(xq^{3/2}x/p)$ corresponds to the weight $(-p/(\sqrt{q}x);q)_\infty \tilde{w}_{1/2}$ when $c=3/2$ and hence has the same moments as
the weight $\rho_0$ in \eqref{chihara-weight}.  

\subsection{Little $q$-Laguerre polynomials}
The monic little $q$-Laguerre polynomials are given by
\[  P_n(x|q) = (-1)^n q^{\binom{n}{2}} (q^{\alpha+1};q)_n\ {}_2\phi_1 \left( \begin{array}{c} q^{-n}, 0 \\ q^{\alpha+1} \end{array} ; q,qx \right), \]
and their orthogonality relation is
\[  \sum_{k=0}^\infty P_n(q^k|q) P_m(q^k|q) q^k \frac{q^{k\alpha}}{(q;q)_k} = 0, \qquad m\neq n, \]
(\cite[\S 14.20]{KoekoekLeskySwarttouw}). This can also be written as
\[   \int_0^1 P_n(x|q)P_m(x|q) x^\alpha (qx;q)_\infty \, d_qx = 0, \qquad m\neq n, \]
where Jackson's $q$-integral is 
\begin{equation} \label{q-Jackson}
  \int_0^1 f(x)\, d_qx = (1-q) \sum_{k=0}^\infty q^k f(q^k).  
\end{equation}
We can therefore say that the little $q$-Laguerre polynomials are orthogonal on the $q$-lattice $\{q^k \ | \ k=0,1,2,\ldots\}$ for the
weight function $w(x) = x^\alpha (qx;q)_\infty$. The little $q$-Laguerre polynomials are intimately related to the $q$-Laguerre polynomials:
\[   L_n^{(\alpha)}(x;q) = (-1)^n P_n(-x|q^{-1}), \]
and from this we easily find the recurrence coefficients from those of the $q$-Laguerre polynomials by changing $q$ to $1/q$ and
changing the sign in $b_n$:
\begin{align*}
b_n &= q^{n} \bigl( 1+q^\alpha - q^{n+\alpha}(1+q) \bigr), \\
a_n^2 &= q^{2n+\alpha-1} (1-q^n)(1-q^{n+\alpha}) .
\end{align*}
The moment problem for the little $q$-Laguerre polynomials is determinate. 

\subsection{Semiclassical extensions: main results} \label{section:main}

In this paper we will consider semiclassical extensions of the
Stieltjes-Wigert, $q$-Laguerre and little $q$-Laguerre weights by replacing the variable $x$ by $x^2$ and the parameter $q$ by $q^2$, but still restricting the weight to the positive real axis. This is very similar to the first semiclassical extension of the Hermite polynomials
(with weight $w(x) = e^{-x^2}$ on $\mathbb{R}$) to orthogonal polynomials for the Freud weight $w(x) = e^{-x^4}$ on $\mathbb{R}$, which leads to
recurrence coefficients satisfying the discrete Painlev\'e I equation (d-$\textrm{P}_{\textrm{I}}$), see e.g., \cite{VanAssche} or the semiclassical extension of the Laguerre polynomials (with weight $w(x)=x^\alpha e^{-x}$ on $[0,\infty))$ to $w(x)=x^\alpha e^{-x^2+tx}$ on $[0,\infty)$, which leads to an asymmetric discrete Painlev\'e IV equation
(d-$\textrm{P}_{\textrm{IV}}$), see \cite{BoelenVA}.
Our main interest is in the recurrence coefficients of the corresponding orthogonal polynomials. We will show that they are related to solutions
of the $q$-discrete Painlev\'e III ($q$-$\textrm{P}_{\textrm{III}}$) and Painlev\'e V ($q$-$\textrm{P}_{\textrm{V}}$) equations.

After this paper was submitted, we learned that Ormerod  considered simple generalizations of the Big $q$-Laguerre polynomials  
and their relation with $q$-$\textrm{P}_{\textrm{V}}$ in \cite[\S 5]{Ormerod}. His results suggest that it would be of interest to investigate semi-classical extensions
of the Big $q$-Laguerre polynomials as well.

In Section \ref{section:Stieltjes-Wigert} we will prove the following result:
\begin{theorem}  \label{thm:1}
The recurrence coefficients of the orthogonal polynomials for the
semiclassical Stieltjes-Wigert weight
$$
w(x) = \frac{x^\alpha}{(-x^2;q^2)_\infty(-q^2/x^2;q^2)_\infty},
\qquad x \in [0,\infty),\  \alpha \in \mathbb{R}
$$
can be found by
\begin{align*}
  q^{2n+\alpha}b_n^2x_n &= x_{n+1}+q^{2n+2\alpha}x_{n-1}\left(x_n+q^{-n-\alpha}\right)^2+2\left(x_n+q^{-\alpha}\right) 
\intertext{and}
  a_n^2 &= q^{1-n}x_n+q^{-2n-\alpha+1},
\end{align*}
where $x_n$ satisfies the $q$-discrete Painlev\'e III equation
\begin{equation}  \label{qPIII-thm1}
x_{n-1}x_{n+1} = \frac{(x_n+q^{-\alpha})^2}{(q^{n+\alpha}x_n+1)^2}.
\end{equation}
\end{theorem}

In Section \ref{section:q-Laguerre} we will prove
\begin{theorem} \label{thm:2}
The recurrence coefficients of the orthogonal polynomials for the
semiclassical $q$-Laguerre weight
$$
w(x)= \frac{x^\alpha
(-p/x^2;q^2)_\infty}{(-x^2;q^2)_\infty(-q^2/x^2;q^2)_\infty},\qquad x\in [0,\infty),\ p\in [0,q^{-\alpha}), \alpha \geq 0
$$
can be found by
\begin{align*}
  b_n^2q^{2n+\alpha}z_n^2 &=
  z_nz_{n+1}-1+q^{2n+2\alpha}\left(\sqrt{pq^{-2-\alpha}}z_n+q^{-n-\alpha}\right)^2\left(z_nz_{n-1}-1\right)\\
 & \qquad + 2\left(z_n+\sqrt{q^{2-\alpha}/p}\right)\left(z_n+\sqrt{pq^{\alpha-2}}\right),
\intertext{and}
  a_n^2 &= q^{-n+1}z_n\sqrt{pq^{-2-\alpha}}+q^{-2n-\alpha+1},
\end{align*}
where $z_n$ satisfies the $q$-discrete Painlev\'e V equation
$$
(z_nz_{n-1} -1)(z_nz_{n+1}-1)=
\frac{\left(z_n+\sqrt{q^{2-\alpha}/p}\right)^2\left(z_n+\sqrt{pq^{\alpha-2}}\right)^2}
{\left(q^{n+\alpha/2-1}\sqrt{p}z_n+1\right)^2}.
$$
\end{theorem}
Note that this result reduces to Theorem \ref{thm:1} if we take $x_n=z_n\sqrt{p}q^{-1-\alpha/2}$ and let $p \to 0$.
The result for $p=q^2$ was already obtained in \cite[Thm.~1.3]{BoelenVA}.

In Section \ref{section:little-q-Laguerre} we will prove
\begin{theorem}  \label{thm:3}
The recurrence coefficients of the orthogonal polynomials  for the weight
\[   w(x) = x^\alpha (q^2x^2;q^2)_\infty, \qquad \alpha >0  \]
on the positive exponential lattice $\{q^k \ | \ k=0,1,2,3,\ldots\}$ are given by
\[      a_n^2 = q^{n+\alpha/2-1}(x_n - q^{n+\alpha/2}), \qquad n \geq 1,  \]
and
\begin{multline*}  b_n^2 x_n^2q^{-2n-\alpha} = 1 - x_nx_{n+1}-q^{-2n}(x_nx_{n-1}-1)(x_nq^{-\alpha/2}-q_n)^2 \\
     -\ 2(x_n-q^{\alpha/2})(x_n-q^{-\alpha/2}), \qquad n \geq 1,
\end{multline*}
with $b_0^2=1-q^{\alpha/2}x_1$, where $x_n$ satisfies the $q$-discrete Painlev\'e V equation
\begin{equation}   \label{qPV-3}
 (x_nx_{n-1}-1)(x_nx_{n+1}-1) = \frac{q^{2n+\alpha}(x_n-q^{\alpha/2})^2(x_n-q^{-\alpha/2})^2}{(x_n-q^{n+\alpha/2})^2}, 
\end{equation}
with initial values $x_0=q^{\alpha/2}$ and $x_1=q^{-\alpha/2}(1-\mu_1^2/\mu_0^2)$, where
$\mu_0$ and $\mu_1$ are the first two moments of the weight $w$.  
\end{theorem}

\section{Ingredients of the proofs}

\subsection{Discrete Painlev\'e equations}
Our Theorems \ref{thm:1}--\ref{thm:3} express the recurrence coefficients as solutions of $q$-discrete Painlev\'e equations.
These equations  are second order, nonlinear difference
equations which have a continuous Painlev\'e equation as a
continuous limit. They pass an integrability test called singularity
confinement \cite{GrammaticosRamaniPapageorgiou}. This integrability detector is the discrete analogue of the
Painlev\'e property for differential equations. 
For an overview of discrete Painlev\'e equations, see
\cite{GrammaticosRamani2}. Three types of discrete Painlev\'e equations are
distinguished: those of additive nature (the independent variable
$n$ entering through $\alpha n+\beta$) are usually denoted by d-P or
$\delta$-P, those of multiplicative nature, where the independent
variable $n$ enters the equation in an exponential way, which are denoted by $q$-P, and an elliptic
equation. Sakai gave a classification of discrete Painlev\'e equations in terms of affine Weyl groups, starting from the
exceptional Weyl group $E_8$ (Sakai \cite{Sakai}, see also \cite[\S 6, p. 293--300]{GrammaticosRamani2}).
In this paper we are dealing with $q$-discrete Painlev\'e III ($q$-$\textrm{P}_{\textrm{III}}$) 
\begin{equation}  \label{qPIII}
     x_{n+1}x_{n-1} = \frac{(x_n+a)(x_n+b)}{(1+cq^nx_n)(1+dq^nx_n)}, 
\end{equation}
and $q$-discrete Painlev\'e V ($q$-$\textrm{P}_{\textrm{V}}$)
\begin{equation}   \label{qPV}
    (x_{n+1}x_n-1)(x_nx_{n-1}-1) = cdq^{2n} \frac{(x_n-a)(x_n-1/a)(x_n-b)(x_n-1/b)}{(x_n-cq^n)(x_n-dq^n)} ,
\end{equation}
\cite[p.~269]{GrammaticosRamani2}, \cite[p.~720]{VanAssche}. We will find special cases of these $q$-discrete Painlev\'e equations with
$a=b$ and $c=d$.

\subsection{Technique of ladder operators}
In what follows we will be interested in retrieving the recurrence
coefficients of semiclassical extensions of the Stieltjes-Wigert polynomials, the $q$-Laguerre polynomials and the little $q$-Laguerre
polynomials. We will use the technique of ladder operators for $q$-orthogonal
polynomials which was introduced by Chen and Ismail in \cite{ChenIsmail}. They consider monic polynomials orthogonal with
respect to an absolutely continuous measure, defined through a
weight function $w$, on the positive real axis. The potential $u$ is
defined as
\begin{equation}\label{potential}
u(x) = -\frac{D_{q^{-1}}w(x)}{w(x)}
\end{equation}
where $D_q$ is the $q$-analogue of the difference operator,
$$
(D_qf)(x) = \begin{cases} 
   \displaystyle \frac{f(x)-f(qx)}{x(1-q)}  &\textrm{if $x \neq 0$},\\
   f'(0)                      & \textrm{if } x=0. 
             \end{cases}
$$
The main result involves two entities
\begin{equation}\label{Aqhalfline}
A_n(x)  =  \gamma_n^2 \int_0^\infty \frac{u(qx)-u(y)}{qx-y}
P_n(y)P_n(y/q)w(y)\, dy
\end{equation}
\begin{equation}\label{Bqhalfline}
B_n(x) = \gamma_{n-1}^2 \int_0^\infty
\frac{u(qx)-u(y)}{qx-y}P_n(y)P_{n-1}(y/q)w(y)\, dy ,
\end{equation}
where
\[   \frac{1}{\gamma_n^2} = \int_0^\infty P_n^2(x)w(x)\, dx, \]
which appear in the lowering relation \cite[Theorem 1.1]{ChenIsmail}
$$
D_qP_n(x) = a_n^2 A_nP_{n-1}(x) - B_nP_n(x).
$$
We can make use of the following identities \cite[Eqs. (1.15)--(1.16)]{ChenIsmail}
\begin{equation}
\int_0^\infty u(y) P_n(y)P_n(y/q)w(y)\, dy = 0
\end{equation}
and
\begin{equation}
\gamma_n^2 \int_0^\infty u(y)P_{n+1}(y)P_n(y/q)w(y)\, dy =
\frac{1-q^{n+1}}{1-q}q.
\end{equation}
Furthermore, the quantities $A_n$ and $B_n$ satisfy two compatibility relations
\begin{equation}\label{rel1}
B_{n+1}(x) +B_n(x) = (x-b_n)A_n(x) + x(q-1)\sum_{j=0}^nA_j(x)-
u(qx),
\end{equation}
\begin{equation}\label{rel2}
1+(x-b_n)B_{n+1}(x) - (qx-b_n)B_n(x) =
a_{n+1}^2A_{n+1}(x) - a_n^2A_{n-1}(x).
\end{equation}
These equations will enable us to find expressions for the
recurrence coefficients $b_n, a_n^2$.

Ismail \cite{Ismail} shows that the ladder operators also work for orthogonal polynomials on the exponential lattice $\{q^n\, | \, n =0,1,2,\ldots\}$,
in particular for the $q$-Jackson integral \eqref{q-Jackson}. He shows that the same formulas hold for the expressions 
\begin{equation}\label{Aq^n}
A_n(x)  =  \gamma_n^2 \int_0^1 \frac{u(qx)-u(y)}{qx-y}
P_n(y)P_n(y/q)w(y)\, d_qy
\end{equation}
and
\begin{equation}\label{Bq^n}
B_n(x) = \gamma_{n-1}^2 \int_0^1
\frac{u(qx)-u(y)}{qx-y}P_n(y)P_{n-1}(y/q)w(y)\, d_qy .
\end{equation}

\section{Semiclassical Stieltjes-Wigert weight} \label{section:Stieltjes-Wigert}

In this section we will prove Theorem \ref{thm:1} using the compatibility relations for the ladder operators.
For the weight
\[   w(x)= \frac{x^\alpha}{(-x^2;q^2)_\infty (-q^2/x^2;q^2)_\infty}, \qquad x \in [0,\infty), \]
we have
\begin{equation} \label{Pearson1}
   w(x/q)=q^{2-\alpha} w(x)/x^2,  
\end{equation}
and hence the potential, defined in \eqref{potential}, is given by
\[   u(x) = \frac{q}{1-q} \left( \frac{1}{x} - \frac{q^{2-\alpha}}{x^2} \right), \]
so that
\[  \frac{u(qx)-u(y)}{qx-y} = -\frac{u(y)}{qx} + \frac{q^{1-\alpha}}{(1-q)x^2y^2} + \frac{q^{-\alpha}}{(1-q)x^3y} . \] 
The formula \eqref{Aqhalfline} then becomes
\[   A_n(x) = \frac{T_n}{(1-q)x^3} + \frac{\gamma_n^2 q^{1-\alpha}}{(1-q)x^2} \int_0^\infty P_n(y)P_n(y/q) \frac{w(y)}{y^2}\, dy, \]
with
\[  T_n = \gamma_n^2 q^{-\alpha} \int_0^\infty P_n(y)P_n(y/q) \frac{w(y)}{y}\, dy , \]
The relation \eqref{Pearson1} gives 
\begin{align*}
   \int_0^\infty P_n(y)P_n(y/q) \frac{w(y)}{y^2}\, dy &= q^{\alpha-2} \int_0^\infty P_n(y)P_n(y/q)w(y/q)\, dy  \\
             &=  q^{\alpha-1} \int_0^\infty P_n(qx)P_n(x)w(x)\, dx \\
             &= \frac{q^{n+\alpha-1}}{\gamma_n^2}, 
\end{align*}
where we have used $P_n(qx) = q^nP_n(x) + \cdots$.  Hence we have that $A_n$ is the rational function
\begin{equation}  \label{Aq1}
   A_n(x) = \frac{T_n}{(1-q)x^3} + \frac{q^n}{(1-q)x^2}.  
\end{equation}
If we write $P_n(x) = x^n + \delta_n x^{n-1} + \cdots$, then \eqref{Pearson1} and a calculation similar as before give
\begin{equation}  \label{T}
   T_n = q^{n-1}(\delta_n - q \delta_{n+1}).  
\end{equation}
 From this we easily find
\[   \sum_{j=0}^n T_j = -q^n\delta_{n+1}.  \]
The formula \eqref{Bqhalfline} gives the rational function
\begin{equation}  \label{Bq1}
   B_n(x) = - \frac{1-q^n}{(1-q)x} + \frac{r_n}{(1-q)x^2} + \frac{t_n}{(1-q)x^3},
\end{equation}
with
\begin{align*}
   r_n &= q^{1-\alpha} \gamma_{n-1}^2 \int_0^\infty P_n(y)P_{n-1}(y/q) \frac{w(y)}{y^2}\, dy, \\
   t_n &= q^{-\alpha} \gamma_{n-1}^2  \int_0^\infty P_n(y)P_{n-1}(y/q) \frac{w(y)}{y}\, dy, 
\end{align*}
where $r_0=t_0=0$. Using $P_n(y)=y^n + \delta_ny^{n-1} + \cdots$ and \eqref{Pearson1} one finds
\begin{equation}  \label{r}
  r_n = q^{n-1}(1-q) \delta_n.  
\end{equation}
Inserting \eqref{Aq1} and \eqref{Bq1} into the compatibility relation \eqref{rel1} and comparing powers of $x$ gives
\begin{subequations}
\begin{align}
   t_{n+1}+t_n &= -b_nT_n + q^{-\alpha},  \label{tT}  \\
   r_{n+1}+r_n &= -b_nq^n + T_n + (q-1) \sum_{j=0}^n T_j,  \label{rT}
\end{align}
\end{subequations}
and the second compatibility relation \eqref{rel2} gives
\begin{subequations}
\begin{align}
     -b_nt_{n+1}+b_n t_n &= a_{n+1}^2 T_{n+1} -a_n^2 T_{n-1}, \label{tTT} \\
     t_{n+1}-qt_n-b_nr_{n+1} + b_n r_n &= a_{n+1}^2q^{n+1} -a_n^2 q^{n-1},  \label{tr}  \\
      b_nq^n(1-q) + r_{n+1} -qr_n &= 0.   \label{rr}
\end{align}
\end{subequations}
The equations \eqref{rT} and \eqref{rr} give $b_n = \delta_n-\delta_{n+1}$, but this relation is well known and
follows by comparing the coefficients of $x^{n}$ in the three term recurrence relation \eqref{TTR}. 
When multiplying \eqref{tTT} by $T_n$  and using \eqref{tT}, one finds
\[     (t_{n+1}+t_n-q^{-\alpha})(t_{n+1}-t_n) = a_{n+1}^2 T_{n+1}T_n - a_n^2 T_{n}T_{n-1}. \]
Summing over $n$ gives
\begin{equation}  \label{TT}
   a_n^2T_nT_{n-1} = t_n(t_n-q^{-\alpha}).
\end{equation}
 From \eqref{r} and \eqref{T} one has
\[   r_{n+1}-r_n = (1-q)(q^n\delta_{n+1} - q^{n-1}\delta_n) = -(1-q)T_n. \]
Inserting this in \eqref{tr} gives
\[  t_{n+1}-qt_n + (1-q)b_nT_n = a_{n+1}^2 q^{n+1} - a_n^2 q^{n-1} . \]
Use \eqref{tT} to obtain
\[   q^{-\alpha}(1-q) - t_n + qt_{n+1} = a_{n+1}^2 q^{n+1}-a_n^{2} q^{n-1}.  \]
Multiply this by $q^n$, then summing over $n$ gives
\begin{equation} \label{at}
  a_n^2 q^{2n-1}  = q^{-\alpha}(1-q^n) + q^n t_n, 
\end{equation}
which expresses $t_n$ in terms of $a_n^2$. We use this to get rid of $t_n$ and $t_{n+1}$ in  \eqref{tTT}:
\[  a_{n+1}^2 (T_{n+1} + q^n b_n) - a_n^2 (T_{n-1}+q^{n-1} b_n)  = q^{-n-1-\alpha} b_n(1-q). \]
Replacing $b_n$ and $T_{n\pm 1}$ by their expressions involving $\delta_n$ gives
\[ b_n(1-q)q^{-\alpha} = q^{2n+1} a_{n+1}^2 (\delta_n - q\delta_{n+2}) - q^{2n-1} a_n^2 (\delta_{n-1}-q\delta_{n+1}). \]
Summing over $n$ gives
\begin{equation}  \label{ba}
 q^{2n+\alpha-1} a_n^2 (\delta_{n-1}-q\delta_{n+1}) = (1-q) \sum_{j=0}^{n-1} b_j = -(1-q) \delta_n. 
\end{equation}
 From \eqref{rT} one finds
\[ q^nb_n = T_n - q^{n-1}(1-q)\delta_n, \]
and this gives for \eqref{ba}
\[  b_n q^n(1-q^{2n+\alpha-1}a_n^2) = T_n + q^{2n+\alpha} b_n T_{n-1}.  \]
We can use this expression in \eqref{tT} to find
\[   (t_{n+1}+t_n-q^{-\alpha})(1-q^{2n+\alpha-1}a_n^2) q^n = -T_n^2-q^{2n+\alpha} a_n^2 T_nT_{n-1}.  \]
The last term in this expression can be replaced using \eqref{TT}, and by using \eqref{at} one has
\[  T_n^2 = q^{2n+\alpha}(t_n-q^{-\alpha})(t_{n+1}-q^{-\alpha}).  \]
We can use this in \eqref{TT} and together with \eqref{at} we find an expression in $t_n$, $t_{n-1}$ and $t_{n+1}$ only
\[  (q^{-\alpha}(1-q^n)+q^nt_n)^2q^{2\alpha}(t_{n+1}-q^{-\alpha})(t_{n-1}-q^{-\alpha}) (t_n-q^{-\alpha})^2 = t_n(t_n-q^{-\alpha}), \]
where we can divide both sides by $(t_n-q^{-\alpha})^2$ to obtain
\[ t_n^2 = q^{2\alpha}\bigl(q^{-\alpha}(1-q^n)+q^nt_n\bigr)^2(t_{n+1}-q^{-\alpha})(t_{n-1}-q^{-\alpha}). \]
If we define $x_n = t_n-q^{-\alpha}$ then we get \eqref{qPIII-thm1} from Theorem \ref{thm:1}. This is $q$-$\textrm{P}_{\textrm{III}}$ 
\eqref{qPIII} with $a=b=q^{-\alpha}$ and $c=d=q^{\alpha}$. The recurrence coefficients are expressed in terms of $x_n$ in the following way:
\begin{align*}
    q^{2n-1}a_n^2 &= q^nx_n+q^{-\alpha}, \qquad n \geq 0 \\
    q^{2n+\alpha} b_n^2 x_n &= x_{n+1} + q^{2n+2\alpha}x_{n-1}(x_n+q^{-n-\alpha})^2 + 2(x_n+q^{-\alpha}), \qquad n \geq 1, 
\end{align*}
where $x_0=-q^{-\alpha}$ and $x_1= b_0^2$ and $b_n$ is to be taken positive.   This proves Theorem \ref{thm:1}.

\section{Semiclassical $q$-Laguerre weight} \label{section:q-Laguerre} 
In this section we will look at a semiclassical variation on the
weight functions \eqref{askey-weight} and \eqref{chihara-weight}, namely
$$
w(x)= \frac{x^\alpha(-p/x^2;q)_\infty}{(-x^2;q^2)_\infty(-q^2/x^2;q^2)_\infty},
\qquad x\in [0,\infty), p\in [0,q^{-\alpha}), \alpha \geq 0.
$$
The case $p=q^2$ is a semiclassical generalization of the
$q$-Laguerre weight $x^\alpha /(-x;q)_\infty$ described in \S \ref{sec:1.2} and \cite{BoelenVA}.

We denote the sequence of monic orthogonal polynomials by
$(P_n)_{n\geq 0}$. When we compute $w(x/q)$ then we find
\begin{equation}\label{breuk}
\frac{w(x)}{p+x^2}= \frac{w(x/q)}{q^{2-\alpha}}.
\end{equation}
The potential for this weight $w$ is
$$
u(x) = \frac{q}{1-q}\left(\frac{1}{x}-\frac{q^{2-\alpha}}{x(p+x^2)}
\right),
$$
from which we find
$$
\frac{u(qx)-u(y)}{qx-y} =
-\frac{1}{qx}u(y)+\frac{q^{3-\alpha}}{1-q}\frac{1}{(p+q^2x^2)(p+y^2)}+\frac{q^{2-\alpha}}{1-q}\frac{y}{x(p+q^2x^2)(p+y^2)}.
$$
We now look for expressions for the functions $A_n$ and $B_n$ of the lowering relation. For $A_n$ we have
\begin{equation}
A_n(x) =
\frac{q^2}{1-q}\frac{T_n}{x(p+q^2x^2)}+\frac{q^{n+2}}{1-q}\frac{1}{p+q^2x^2}, 
\end{equation}
with
$$
T_n =  q^{-\alpha} \gamma_n^2 \int_0^\infty
yP_n(y)P_n(y/q)\frac{w(y)}{p+y^2}\, dy ,
$$
while
\begin{equation}
B_n(x) =
-\frac{1-q^n}{1-q}\frac{1}{x}+\frac{q^2}{1-q}\frac{r_n}{p+q^2x^2}+\frac{q^2}{1-q}\frac{t_n}{x(p+x^2)}, 
\end{equation}
with
$$
r_n = q^{1-\alpha} \gamma_{n-1}^2 \int_0^\infty
P_n(y)P_{n-1}(y/q)\frac{w(y)}{p+y^2}\, dy
$$
and
$$
t_n = q^{-\alpha}\gamma_{n-1}^2 \int_0^\infty
yP_n(y)P_{n-1}(y/q)\frac{w(y)}{p+y^2}\, dy.
$$
It is clear that $r_0=t_0=0$ as $P_{-1}=0$. Bearing in mind that we
write the polynomials $P_n$ as
$$
P_n(x) = x^n + \delta_n x^{n-1} +  \cdots,
$$
it's easy to see that
\begin{equation}  \label{T}
T_n = q^{n-1}\left(\delta_n-q \delta_{n+1} \right)
\end{equation}
and
\begin{equation}\label{r}
r_n = q^{n-1}(1-q)\delta_n.
\end{equation}
We used \eqref{breuk} in these computations and the substitution
$u=y/q$ in the integrals and we wrote $uP_n(qu)$ and $P_n(qu)$ in
their Fourier expansion to obtain these results. Using the newly found
expression for $T_n$, we can easily compute the following sum which
will appear in the compatibility relations:
$$
\sum_{j=0}^n T_j = -q^n\delta_{n+1}.
$$
Multiplying the first compatibility relation \eqref{rel1} by
$q^{-2}(1-q)x(p+x^2$) and comparing coefficients of powers of $x$,
we get two relations
\begin{subequations}
\begin{align}\label{10}
(t_{n+1}+t_n)+pq^{-2}(-1+q^n+q^{n+1}) &= -b_n T_n +q^{-\alpha}, \\
r_{n+1}+r_n &= -b_nq^n+T_n + (q-1)\sum_{j=0}^nT_j.  \label{11}
\end{align}
\end{subequations}

From the second compatibility relation \eqref{rel2} we get in the
same way the following three relations:
\begin{subequations}
\begin{align}\label{20}
-b_n(t_{n+1}-t_n)+b_nq^n(1-q)pq^{-2} &=
a_{n+1}^2T_{n+1}-a_n^2T_{n-1}, \\
 \label{21}
t_{n+1}-qt_n -b_nr_{n+1}+b_nr_n &=
a_{n+1}^2q^{n+1}-a_n^2q^{n-1}, \\
 \label{22}
b_n(1-q^{n+1})+r_{n+1}-b_n(1-q^n)-qr_n &=0.
\end{align}
\end{subequations}
The equations \eqref{11} and \eqref{22} are trivial in the light of
\eqref{T} and \eqref{r} as
$$
b_n = \delta_n-\delta_{n+1}.
$$
We elaborate on \eqref{21} by noticing that
\begin{align*}
r_{n+1}-r_n &= (1-q)\left( q^n\delta_{n+1}-q^{n-1}\delta_n \right)\\
     &=  -(1-q)T_n.
\end{align*}
Inserting this in \eqref{21} we obtain
$$
t_{n+1}-qt_n + (1-q)b_nT_n = a_{n+1}^2q^{n+1}-a_n^2q^{n-1}.
$$
We use \eqref{10} to eliminate $b_nT_n$ and we obtain
$$
q^{-\alpha}(1-q) - t_n +qt_{n+1}+(1-q)pq^{-2}(1-q^n-q^{n+1}) =
a_{n+1}^2q^{n+1}-a_n^2q^{n-1}.
$$
This equation has the integrating factor $q^n$:
$$
q^{n-\alpha}(1-q)+q^{n+1}t_{n+1}-q^t_n
+(1-q)pq^{-2}(q^n-q^{2n}(1+q))= a_{n+1}^2q^{2n+1}-a_n^2q^{2n-1}.
$$
Summing over $n$ (a telescopic sum) we get
\begin{equation}\label{21'}
a_n^2q^{2n-1}=\left(q^{-\alpha}+pq^{-2}\right)(1-q^n)-pq^{-2}(1-q^{2n})+q^nt_n,
\end{equation}
which expresses $t_n$ in terms of $a_n^2$.

The right hand side of \eqref{20} has the integrating factor $T_n$. When
multiplying this identity by $T_n$, the factor $b_nT_n$
appears twice on the left hand side. We use \eqref{10} to eliminate
this factor when it accompanies the term $t_{n+1}-t_n$ and \eqref{21}
when it accompanies the term $pq^n$. Then we get
\begin{align*}
a_{n+1}^2T_{n+1}T_n -a_n^2T_nT_{n-1} &= 
\Bigl(t_{n+1}+t_n-q^{-\alpha}+p(q^{-2}-1+q^n+q^{n+1})\Bigr) (t_{n+1}-t_n)\\
 & \quad  +\ pq^{-2}q^n\Bigl( a_{n+1}^2q^{n+1}-a_n^2q^{n-1}-t_{n+1}+qt_n \Bigr).
\end{align*}
Using \eqref{21'} in the factor of $pq^{-2+n}$ we get
\begin{align*}
a_{n+1}^2T_{n+1}T_n -a_n^2T_nt_{n-1} &= t_{n+1}^2-t_n^2 -
(t_{n+1}-t_n)(q^{-\alpha}+pq^{-2})+ 2pq^{-2}q^{n+1}t_{n+1} \\
 & \quad  -\ 2pq^{-2}q^nt_n + pq^{-2} (q^{-\alpha}+pq^{-2})(1-q)q^n-p^2q^{-4}
q^{2n}(1-q^2).
\end{align*}
Summing over $n$ (telescopic sum) and taking into account that
$t_0=a_0^2=0$, we get
\begin{align}\label{20'}
  a_n^2T_nT_{n-1} &= t_n(t_n-q^{-\alpha}-pq^{-2}+2pq^{-2}q^n)   \nonumber \\
& \quad   +\ pq^{-2}(q^{-\alpha}+pq^{-2})(1-q^n)-p^2q^{-4}(1-q^{2n}).
\end{align}
We  now use \eqref{21'} to eliminate $t_n$ and $t_{n+1}$ in \eqref{20}
and obtain
$$
a_{n+1}^2 ( T_{n+1} +q^n b_n ) -a_n^2 ( T_{n-1}+q^{n-1}b_n ) = q^{-n-1}b_n(1-q).
$$
Replacing $b_n$ and $T_n$ in the left hand side by their
expressions in $\delta_n$, we get
$$
b_n(1-q)q^{-\alpha} =
q^{2n+1}a_{n+1}^2(\delta_n-q\delta_{n+2})-q^{2n-1}a_n^2(\delta_{n-1}-q\delta_{n+1}).
$$
Taking a telescopic sum, we obtain
\begin{equation}\label{p_1}
q^{2n+\alpha-1}a_n^2\left(\delta_{n-1} -q\delta_{n+1} \right) =
(1-q)\sum_{j=0}^{n-1} b_j = -(1-q)\delta_n. 
\end{equation}
We can now look at $b_n$:
$$
q^nb_n = T_n -q^{n-1}(1-q)\delta_n.
$$
Replacing $\delta_n$ using \eqref{p_1}, we get
\begin{align*}
q^n b_n &=  T_n + q^{3n+\alpha-2}a_n^2 ( \delta_{n-1}-q\delta_{n+1} )  \\
 &=  T_n + q^{3n+\alpha-1}a_n^2b_n+q^{2n+\alpha}a_n^2T_{n-1}.
\end{align*}
Collecting all terms with $b_n$ then gives
\begin{equation}\label{alpha}
b_nq^n ( 1-q^{2n+\alpha-1}a_n^2 )=T_n +q^{2n+\alpha}a_n^2 T_{n-1}  .
\end{equation}
We can use this expression in \eqref{10}
$$
\Bigl( t_{n+1}+t_n-q^{-\alpha}+pq^{-2}(-1+q^n+q^{n+1}) \Bigr) (1-q^{2n+\alpha-1}a_n^2)q^n
= -T_n^2-q^{2n+\alpha}a_n^2T_nT_{n-1}.
$$
The last term on the right hand side can be replaced using
\eqref{20'}. We then get an expression for $T_n^2$ in terms of $t_n$:
\begin{align*}    \label{beta}
T_n^2 &= -q^{2n +\alpha}\Big( t_n(q^{-\alpha}+pq^{-2}-pq^{-2}q^{n+1})+t_{n+1}(-tn+q^{-\alpha} +
pq^{-2}-pq^{-2}q^n)\nonumber\\ 
 &\quad   +\ pq^{-2}q^{n-\alpha+1}-p^2q^{-4}(1-q^n-q^{n+1}+q^{2n+1}) 
  -q^{-\alpha}(q^{-\alpha}+pq^{-2}-pq^{-2}q^n) \Big) .
\end{align*}
This expression can be simplified by using the substitution $y_n = t_n
-q^{-\alpha}-pq^{-2}+pq^{-2}q^n$, giving
\begin{equation}\label{T_n^2}
T_n^2 = q^{2n+\alpha}(y_ny_{n+1}-pq^{-2}q^{-\alpha}).
\end{equation}
We can use this expression in \eqref{20'}, which after squaring becomes
$$
 a_n^4T_n^2 T_{n-1}^2 = \left( y_n +q^{-\alpha} \right)^2 \left(
y_n +pq^{-2} \right)^2 ,
$$
and together with \eqref{21'} we are left with an expression in
$y_{n-1},y_n$ and $y_{n+1}$ only:
$$
(y_ny_{n-1}-pq^{-2}q^{-\alpha})(y_ny_{n+1}-pq^{-2}q^{-\alpha}) =
\frac{(y_n+q^{-\alpha})^2 (y_n+pq^{-2})^2}{(q^{n+\alpha}y_n+1)^2}.
$$
If we now use the substitution $z_n = y_n/\sqrt{pq^{-2-\alpha}}$
(assuming $p \neq 0$), we get
$$
(z_nz_{n-1} -1)(z_nz_{n+1}-1)=
\frac{\left(z_n+\sqrt{q^{2-\alpha}/p}\right)^2\left(z_n+\sqrt{pq^{\alpha-2}}\right)^2}
{\left(q^{n+\alpha/2-1}\sqrt{p}z_n+1\right)^2},
$$
which is an instance of $q$-P$_{\rm V}$  \eqref{qPV}
with $a=b=c=d=-\sqrt{q^{2-\alpha}/p}$ but with $q$ replaced by $1/q$. 
The initial conditions are
$$
z_0 =-\sqrt{q^{2-\alpha}/p}, \hspace{3mm} z_1 =
\frac{\mu_2\mu_0-\mu_1^2-\mu_0^2q^{-\alpha-1}}{\mu_0^2\sqrt{pq^{-\alpha-2}}},
$$
where $\mu_0$ and $\mu_1$ are the first two moments of the weight $w$. Using
\eqref{21'}, \eqref{alpha} and \eqref{beta}, we can express the
recurrence coefficients $b_n, a_n^2$ in terms of $z_n$:
\begin{multline*}
b_n^2q^{2n+\alpha}z_n^2  = 
z_nz_{n+1}-1+q^{2n+2\alpha}\left(\sqrt{pq^{-2-\alpha}}z_n+q^{-n-\alpha}\right)^2\left(z_nz_{n-1}-1\right)\\
+\ 2\left(z_n+\sqrt{q^{2-\alpha}/p}\right)\left(z_n+\sqrt{pq^{\alpha-2}}\right),
\end{multline*}
and
$$
a_n^2 =q^{-n+1}z_n\sqrt{pq^{-2-\alpha}}+q^{-2n-\alpha+1}.
$$
This concludes the proof of Theorem \ref{thm:2}.

\section{Semiclassical little $q$-Laguerre polynomials} \label{section:little-q-Laguerre}

We consider the weight
\[   w(x) = x^\alpha (q^2x^2;q^2)_\infty, \qquad \alpha > 0, \]
on the exponential lattice $\{q^n\, | \, n=0,1,2,\ldots\}$ and consider monic orthogonal polynomials satisfying
\[  \int_0^1 P_n(x)P_m(x) w(x)\, d_qx = \delta_{n,m}/\gamma_n^2, \]
for the $q$-Jackson integral \eqref{q-Jackson}. The potential $u$ satisfies
\[    u(qx) = - \frac{D_qw(x)}{w(qx)} = \frac{1}{1-q} \left( \frac{1-q^{-\alpha}}{x} + q^{-\alpha+2}x \right), \]
which leads to 
\[   \frac{u(qx)-u(y)}{qx-y} = \frac{1}{1-q} \left( q^{-\alpha+1} - \frac{1-q^{-\alpha}}{xy} \right). \]
Observe that $w(0)=w(1/q)=0$ which implies that all the boundary terms in the expressions for $A_n$ and $B_n$ in \cite{Ismail} vanish.
The function $A_n$ in \eqref{Aq^n} is
\begin{equation*} 
 A_n(x) = \gamma_n^2 \frac{q^{-\alpha}-1}{(1-q)x} \int_0^1 P_n(y)P_n(y/q) \frac{w(y)}{y}\, d_qy 
    + \gamma_n^2 \frac{q^{1-\alpha}}{1-q} \int_0^1 P_n(y)P_n(y/q) w(y)\, d_qy . 
\end{equation*}
It is easy to see that the last integral is $q^{-n}/\gamma_n^2$. We then have
\begin{equation}  \label{Aq3}
   A_n(x) = \frac{R_n}{(1-q)x} + \frac{q^{-n-\alpha+1}}{1-q}, 
\end{equation}
where
\[  R_n = \gamma_n^2 (q^{-\alpha}-1) \int_0^1 P_n(y)P_n(y/q) \frac{w(y)}{y}\, d_qy.  \]
In a similar we find that $B_n$ in \eqref{Bq^n} is
\begin{equation}  \label{Bq3}
   B_n(x) = \frac{r_n}{(1-q)x}, 
\end{equation} 
with
\[  r_n = \gamma_{n-1}^2 (q^{-\alpha}-1) \int_0^1 P_n(y)P_{n-1}(y/q) \frac{w(y)}{y}\, d_qy . \]
If we insert \eqref{Aq1}--\eqref{Bq1} into the compatibility relations \eqref{rel1}--\eqref{rel2} and equate the powers of $x$, then
\begin{subequations}
\begin{align}   
    r_{n+1}+r_n &= - b_n R_n - (1-q^{-\alpha}), \label{rR} \\
    b_nq^{-n-\alpha+1} &= R_n + (q-1) \sum_{j=0}^n R_j, \label{bR}
\end{align}
\end{subequations}
and
\begin{subequations}
\begin{align} 
   a_{n+1}^2 R_{n+1} - a_n^2 R_{n-1}  &= -b_n (r_{n+1}-r_n), \label{aRbr} \\
   q^{-n-\alpha} a_{n+1}^2 - q^{-n-\alpha+2}a_n^2 &= r_{n+1}-qr_n+1-q .   \label{ar}
\end{align}
\end{subequations}
We will use these four equations to find expressions for the recurrence coefficients $a_n^2,b_n$ of the orthogonal polynomials.
Multiplying \eqref{ar} by $q^{-n-2}$ and summing over $n$ gives
\begin{equation}  \label{a2r}
    a_n^2 = q^{n+\alpha-1} (r_n+1-q^n),
\end{equation}
with $a_0=r_0=0$. Multiplying \eqref{aRbr} by $R_n$, eliminating $-b_nR_n$ using \eqref{rR}, and summing over $n$ gives
\begin{equation}  \label{aRR}
   a_n^2R_nR_{n-1} = r_n(r_n+1-q^{-\alpha}).
\end{equation}
Next, we use \eqref{a2r} in \eqref{aRbr} and replace $b_n$ by using \eqref{bR} to find
\begin{multline*}
  (1-q^{-n-1})(q^{n+2}R_{n+1}+q^{n+1}R_n) - (1-q^{-n})(q^{n+1}R_n+q^nR_{n-1}) + R_n(1-q) \\
  = r_{n+1} \left(qR_{n+1} + R_n - (1-q) \sum_{j=0}^n R_j \right)
     - r_n \left( qR_n + R_{n-1} - (1-q) \sum_{j=0}^{n-1} R_j \right). 
\end{multline*}
Summing over $n$ (telescoping sum) gives
\[  r_n \left( qR_n + R_{n-1} - (1-q) \sum_{j=0}^{n-1} R_j \right) = (1-q) \sum_{j=0}^{n-1} + (1-q^{-n})(q^{n+1}R_n+q^nR_{n-1}). \]
Now use \eqref{bR} to get
\begin{equation}  \label{bRR}
  b_n(1+r_n)q^{-\alpha-n+1} = q^{n+1}R_n - R_{n-1}(r_n+1-q^n). 
\end{equation}
We use this expression for $b_n$ in \eqref{rR}
\[ q^{-\alpha-n+1} (1+r_n)(r_{n+1}+r_n+1-q^{-\alpha}) = -q^{n+1}R_n^2 + R_nR_{n-1}(r_n+1-q^n). \]
Combining \eqref{a2r} and \eqref{aRR}  we get
\begin{equation}  \label{R^2}
   R_n^2 = -q^{-2n-\alpha}\Bigl( (r_n+1)(r_{n+1}+1)-q^{-\alpha} \Bigr).
\end{equation}
We can use this expression for $R_n^2$ together with \eqref{a2r} in \eqref{aRR} to find
\[  r_n^2(r_n+1-q^{-\alpha})^2 = q^{-2n}(r_n+1-q^n)^2 \Bigl( (r_n+1)(r_{n+1}+1)-q^{-\alpha} \Bigr)
           \Bigl( (r_n+1)(r_{n-1}+1)-q^{-\alpha} \Bigr). \]
Now define $x_n = q^{\alpha/2}(r_n+1)$, then we find \eqref{qPV-3}, which is $q$-$\textrm{P}_{\textrm{V}}$ \eqref{qPV} 
with $a=b=c=d=q^{\alpha/2}$. The initial values are $x_0=q^{\alpha/2}$ and $x_1 = q^{-\alpha/2} (1-b_0^2)$, where
$b_0=\mu_1/\mu_0$ is the ratio of the first two moments of $w$. Using \eqref{a2r} we have
\[   a_n^2 = q^{n+\alpha-1} (x_nq^{-\alpha/2}-q^n), \]
and $b_n$ can be found using \eqref{aRR}--\eqref{R^2} by
\[   b_n^2x_n^2q^{-2n-\alpha} = 1-x_nx_{n+1}-q^{-2n}(x_nx_{n-1}-1)(x_nq^{-\alpha/2}-q^n)^2
 -2(x_n-q^{\alpha/2})(x_n-q^{-\alpha/2}). \]
This proves Theorem \ref{thm:3}.

\section{Acknowledgements}
This research is part of the PhD dissertation \cite{Boelen} of Lies Boelen. We acknowledge support of the Belgian Interuniversity Attraction Poles Programme P7/18, FWO research grant G.0934.13, and KU Leuven research grant OT/12/073.

\begin{verbatim}
Lies Boelen (l.boelen@imperial.ac.uk)
Walter Van Assche (walter@wis.kuleuven.be)
Department of Mathematics
KU Leuven
Celestijnenlaan 200B box 2400
BE-3001 Leuven, BELGIUM
\end{verbatim}

\end{document}